\documentclass[11pt,a4paper]{article}
\usepackage{enumitem}

\usepackage{setspace}

\usepackage{amsmath}
\usepackage{amssymb}
\usepackage{amsfonts}
\usepackage{amsthm}

\usepackage{array}

\usepackage[numbers]{natbib}
\usepackage[page,toc]{appendix}

\usepackage[T1]{fontenc}
\usepackage{lmodern}

\usepackage{upquote}

\usepackage{listings}

\usepackage{algpseudocode}

\usepackage{tikz,pgfplots}
  \usetikzlibrary{fit,arrows.meta}
  \pgfplotsset{compat=1.3}

\usepackage{hyperref}

\newcommand{\FF}{\mathbb{F}}

\newcommand{\QQ}{\mathbb{Q}}

\newcommand{\ZZ}{\mathbb{Z}}

\newcommand{\cA}{\mathcal{A}}

\newcommand{\cO}{\mathcal{O}}
\newcommand{\cP}{\mathcal{P}}

\newcommand{\cPst}{\cP^*}
\newcommand{\cPstres}{\cPst_{\text{res}}}
\newcommand{\cPstpr}{\cP^{*\prime}}
\newcommand{\cPstrespr}{\cPstpr_{\text{res}}}

\newcommand{\Gal}{\operatorname{Gal}}

\newcommand{\braces}[1]{\left\lbrace#1\right\rbrace}

\newcommand{\ceiling}[1]{\left\lceil#1\right\rceil}
\newcommand{\floor}[1]{\left\lfloor#1\right\rfloor}
\newcommand{\abs}[1]{\left\lvert#1\right\rvert}

\def\comment{}
\def\endcomment{}
\long\def\comment#1\endcomment{}

\theoremstyle{plain}

\theoremstyle{definition}

\newtheorem{example*}[lemma]{Example} 
  {%
   \pushQED{\qed}\begin{example*}}
  {\popQED\end{example*}}
\theoremstyle{remark}

\makeatletter
\renewcommand*{\verbatim@font}{\ttfamily\fontseries{m}\selectfont}
\makeatother

\lstdefinelanguage{Magma}{
  morekeywords={end,function,intrinsic,procedure,for,while,repeat,until,do,in,if,else,elif,then,error,assert,require,when,where,is,print,printf,vprint,vprintf,time,declare,verbose,type,attributes,return,continue,break,delete,loop},
  morekeywords=[2]{eq,ne,le,lt,ge,gt,cmpeq,cmpne,not,and,or,subset,meet,join,diff,sdiff,assigned,eval},
  morekeywords=[3]{sub,ncl,func,proc,ideal,elt},
  morekeywords=[4]{StrPadExact,PadExactElt,FldPadExact,FldPadExactElt,RngUPol_FldPadExact,RngUPolElt_FldPadExact,RngMPol_FldPadExact,RngMPolElt_FldPadExact,SetCart_PadExactElt,Tup_PadExactElt,Val_PadExactElt,Val_FldPadElt,Val_RngUPolElt_FldPad,Val_RngMPolElt_FldPad,RngInt,RngIntElt,SetCart,Tup,FldNum,FldNumElt,FldRat,FldRatElt,FldPad,FldPadElt,Getter,BoolElt,SetCart_PadExact,Tup_PadExact},
  sensitive=true,
  morecomment=[l]{//},
  morecomment=[s]{/*}{*/},
  morestring=[b]",
  literate={->}{\(\rightarrow\) }{1},
}
\begin{document}

\title{On enumerating extensions of \(p\)-adic fields \\ with given invariants}
\author{Christopher Doris \\ University of Bristol \\ \texttt{christopher.doris@bristol.ac.uk}}
\date{March 2018}
\maketitle
\begin{abstract}
We give a brief re-exposition of the theory due to Pauli and Sinclair of ramification polygons of Eisenstein polynomials over \(p\)-adic fields, their associated residual polynomials and an algorithm to produce all extensions for a given ramification polygon. We supplement this with an algorithm to produce all ramification polygons of a given degree, and hence we can produce all totally ramified extensions of a given degree.
\end{abstract}


\section{Introduction}
\label{ext-sec-intro}

We fix a \(p\)-adic field \(K\) and an integer \(n \geq 1\) and consider the problem of enumerating all the extensions \(L/K\) of degree \(n\). In their paper \cite{PS}, Pauli and Sinclair describe some invariants of \(L/K\) and give an algorithm to produce all extensions with a given set of invariants. The algorithm works by producing a set of Eisenstein polynomials generating all extensions for a given invariant, and then determining which pairs of polynomials produce isomorphic extensions so that only one extension per isomorphism class is returned. The finer the invariant used, the smaller the set of polynomials produced, and so easier this pairwise search becomes.

In \S\ref{ext-sec-rampol}--\ref{ext-sec-shrink} we give a brief re-exposition of the invariants they describe, using slightly different notation which we hope is easier to follow. For each invariant, we give its definition, prove that it is an invariant, determine which Eisenstein polynomials \(f(x)\) generate an extension with the given invariant, and give an algorithm to enumerate all possiblities for the invariant.

In \S\ref{ext-sec-rampol} we look at the ramification polygon. This is not studied directly in \cite{PS}, which actually starts with a slightly finer invariant. However the ramification polygon is more widely known, which is why we start with it. In \S\ref{ext-sec-finerampol} we look at the \emph{fine ramification polygon} which is the finer invariant studied in \cite[\S3]{PS}, in which it is called the ``ramification polgon''.

In \S\ref{ext-sec-residues} we look at the \emph{fine ramification polygon with residues} which attaches a residue to each point in the polygon. This is equivalent to a pair of a ramification polygon and the invariant \(\cA\) of \cite[\S4]{PS}, but is notationally simpler: \(\cA\) is expressed as a set of residual polynomials, but this structure is actually mostly irrelevant to studying them.

In \S\ref{ext-sec-constcoeff} we extend this invariant to include some information about the constant coefficient of the Eisenstein \(f(x)\). This is essentially the object \(\cA^*\) of \cite[Eq. 4.2]{PS}, although it is not explicitly identified as an invariant.

In \S\ref{ext-sec-shrink} we consider, as in \cite[\S5]{PS}, transformations to the Eisenstein polynomial \(f(x)\) which preserve the extension \(L/K\), allowing us to reduce the number of such polynomials we need to enumerate in order to find all extensions.

We reiterate that the main results in this article are not new, and are all essentially from \cite{PS}. Exceptions to this are made in the footnotes. We do however use different notation and present the proofs in a different manner.

Many of the ideas here were studied by Monge \cite{Monge} to produce, given a totally ramified extension \(L/K\), a finite set of Eisenstein polynomials each generating \(L/K\), forming an invariant for \(L/K\). The same ideas were also used by Sinclair to count the number of extensions with a given invariant \cite{Sinclair}, the count being the number of extensions within an algebraic closure. Monge has also counted the number of extensions of a given degree up to isomorphism \cite{Monge11} using class field theory to count cyclic extensions, and a group theory argument to count conjugacy classes of subgroups.

The invariants and algorithms described here have been implemented \cite{extensionscode} in the Magma computer algebra system \cite{magma}. In \S\ref{ext-sec-impl} we give a few notes on the implementation. For convenience, our package also provides an implementation of \cite[Lemma 4.4]{Sinclair} to count the number of extensions \(L/K\) in some algebraic closure \(\bar K\) with a given invariant.

\subsection{Notation}

\(K\) is a finite extension of \(\QQ_p\). Fix a uniformizer \(\pi \in K\), and let \(v\) denote the valuation on \(\bar K\) such that \(v(\pi)=1\). We denote by \(\cO_K\) the ring of integers of \(K\), \(\FF_K = \cO_K/(\pi)\) the residue class field, and for \(x \in \cO_K\) we let \(\bar x\) be its residue class.

For \(x,y \in \cO_{\bar K}\), we say \(x \equiv y\) iff \(v(x-y)>0\). For \(x,y \in \bar K\), we say \(x \sim y\) iff \(x=y=0\) or \(v(x-y) > v(x) = v(y)\).

\section{Ramification polygon}
\label{ext-sec-rampol}

\subsection{Definition}
\label{ext-sec-rampol-def}

For a monic Eisenstein polynomial \(f(x) = \sum_{i=0}^n f_i x^i \in K[x]\) of degree \(n\), with a root \(\alpha \in \bar K\), its \emph{ramification polynomial} is
\[r(x) := \alpha^{-n} f(\alpha x + \alpha) \in L[x], \qquad  L := K(\alpha).\]
Note that it is monic and \(r(0)=0\). Expanding out we find
\[r(x) = \sum_{j=1}^n r_j x^j \quad\text{where}\quad r_j = \sum_{i=j}^n \binom{i}{j} f_i \alpha^{i-n} \text{ for } 0<j\leq n.\]

Observing that \(n v(\binom{i}{j} f_i \alpha^{i-n}) \equiv i \mod n\), by the ultrametric property of the valuation we deduce that
\[R_j := n v(r_j) = \min_{i=j}^n n(B(i,j)+F_i-1)+i\]
where \(B(i,j):=v\binom{i}{j}\) and \(F_i := v(f_i)\).

We define the \emph{ramification points of \(f\)} to be
\[\cP = \{(j,R_j) \,:\, 1 \leq j \leq n, R_j < \infty\}\]
and we define the \emph{ramification polygon of \(f\)} to be the lower convex hull of these points, denoted \(P\). That is, it is the Newton polygon of \(r(x)\).

\subsection{Invariant}
\label{ext-sec-rampol-invar}

It is well-known that the ramification polygon \(P\) is actually an invariant of the field \(L/K\). One way to see this is to observe that if \(\alpha'\) is any uniformizer for \(L/K\), then \(\alpha' = x_0 + x_1 \alpha + \ldots + x_{n_1} \alpha^{n-1}\) with \(x_i \in \cO_K\), \(v(x_0)>0\), \(v(x_1)=0\) (since \(\cO_L = \cO_K[\alpha]\)), so for any \(K\)-embedding \(\sigma : L \to \bar K\),
\begin{align*}
\alpha' - \sigma \alpha' &= \sum_{i=0}^n x_i (\alpha^i - \sigma \alpha^i) \\
&= (\alpha - \sigma \alpha)(x_1 + \sum_{i=2}^{n-1} x_i(\alpha^{i-1} + \ldots + \sigma\alpha^{i-1})) \\
&\sim x_1 (\alpha - \sigma \alpha)
\end{align*}
and therefore
\[\frac{\alpha' - \sigma\alpha'}{\alpha'} \sim \frac{\alpha - \sigma\alpha}{\alpha}\]
and in particular these have the same valuation. Since these, over all \(\sigma\), are the roots of the ramification polynomials corresponding to \(\alpha'\) and \(\alpha\), and they have the same valuations, then the Newton polygons of the ramification polynomials are the same.

It is related to the ramification filtration of \(\Gal(L/K)\) if \(L/K\) is Galois (e.g. \cite[Ch. IV]{SerLF}), and more generally to the ramification filtration of the \emph{Galois set} \(\Gamma(L/K) = \{\sigma : L \to \bar K\}\) of \(K\) embeddings of \(L\) (e.g. \cite{Helou}). The vertices of the polygon correspond to subfields of \(L/K\) which themselves correspond to fixed fields of elements of the ramification filtration of \(\Gamma(L/K)\):
\[\Gamma_V := \{\sigma \in \Gamma \,:\, x \in \cO_K \implies v(\sigma x - x) > V\} = \{ \sigma \in \Gamma \,:\, v(\sigma \alpha - \alpha) > V\}.\]

\subsection{Validity}
\label{ext-sec-rampol-valid}

From facts about ramification, or alternatively directly from facts about valuations of binomial coefficients, one can show that the interior vertices of \(P\) are of the form \((p^s,*)\), and that there is a vertex at \((p^{v_p(n)},0)\). Hence we notate a ramification polygon by listing its vertices like so:
\[P = [(p^{s_0},J_0),\ldots,(p^{s_u},J_u=0),(n,0)]\]
where \(s_u = v_p(n)\).

Any polygon of the above form is called a \emph{potential ramification polygon}. Any such polygon arising from an Eisenstein \(f(x) \in K[x]\) is called a \emph{valid ramification polygon (over \(K\))}. We now consider which potential ramification polygons are valid.

Observing that \(P\) is the graph of a function \([1,n] \to \QQ\), we let \(P\) also denote this function.

A potential polygon \(P\) is valid if and only if there is Eisenstein \(f(x)\) with \(v(f_i)=F_i\) such that \(R_{p^{s_t}} = J_t\) for \(0\leq t \leq u\) and such that \(R_j \geq P(j)\) for all \(1 \leq j \leq n\). Automatically we have \(R_j \geq 0\) for all \(j\), which rules the face \([(p^{s_u},0),(n,0)]\) out of consideration. From facts about valuations of binomials, we actually have that if \(p^s < j < p^{s+1}\) then \(R_j \geq R_{p^s} \geq P(p^s) > P(j)\), and hence the only points \((j,R_j)\) lying in \(P\) with \(1 \leq j \leq p^{s_u}\) are of the form \((p^s,R_{p^s})\). We deduce that \(P\) is valid if and only if \(R_{p^s} \geq P(p^s)\) for all \(s\) with equality whenever \(s=s_t\).

Fix some \(s\), then \(R_{p^s} \geq P(p^s)\) if and only if for all \(p^s \leq i \leq n\) we have
\[n(B(i,p^s) + F_i - 1) + i \geq P(p^s)\]
which may be rewritten as
\[F_i \geq \ell_P(i,s) := \ceiling{\frac{P(p^s)-i}{n}} - B(i,p^s) + 1.\]

Suppose this is true for \(s=s_t\), then we have equality \(R_{p^{s_t}} = P(p^s) = J_t =: a_t n + b_t\) with \(1 \leq b_t \leq n\)\footnote{Other authors use \(0 \leq b_t < n\), but we can often avoid special cases with this definition.} if and only if we have equality for \(i = b_t\), i.e. \(p^{s_t} \leq b_t \leq n\) and
\[F_{b_t} = \ell_P(b_t,s_t) = a_t - B(b_t, p^{s_t}) + 1.\]

We deduce that \(P\) is valid if and only if \(p^{s_t} \leq b_t \leq n\) and there are \(F_0,F_1,\ldots,F_n \in \ZZ\) satisfying:
\begin{align*}
& F_0 = 1 \\
0 < i < n \implies& F_i \geq 1 \\
& F_n = 0 \\
p^s \leq i \leq n \implies& F_i \geq \ell_P(i,s) \\
& F_{b_t} = \ell_P(b_t,s_t).
\end{align*}

Each condition is either a lower bound or an equality for some \(F_i\). Hence this system is consistent if and only if the equalities for the same \(F_i\) match, and if the lower bounds are satisfied by the equalities. Hence \(P\) is valid if and only if\footnote{This is essentially \cite[Prop. 3.9]{PS}.}
\begin{align*}
& p^{s_t} \leq b_t \\
b_t = n \implies& \ell_P(n,s_t) = 0 & \text{(Ore 1)}\\
& \ell_P(n,s) \leq 0 & \text{(Ore 2)}\\
b_t < n \implies& \ell_P(b_t,s_t) \geq 1 &\text{(Ore 3)} \\
b_r = b_t < n \implies& \ell_P(b_t,s_t) = \ell_P(b_r,s_r) &\text{(Consistency)} \\
p^s \leq b_t < n \implies& \ell_P(b_t,s_t) \geq \ell_P(b_t,s). &\text{(Bounding)} \\
\end{align*}

Furthermore, if \(P\) is valid, then \(f(x)\) has \(P\) as its ramification polygon if and only if \(F_{b_t} = \ell_P(b_t,s_t)\) for all \(t\) and \(F_i \geq \ell_P(i,s)\) whenever \(p^s \leq i \leq n\).\footnote{This is essentially \cite[Prop. 3.10]{PS}.}

Note that the three ``Ore conditions'' are so-called because they imply
\[\min(n B(b_t, p^{s_t}), n B(n, p^{s_t})) \leq J_t \leq n B(n, p^{s_t})\]
which for \(t=0\) is the bound \(\min(nv(b_0),nv(n)) \leq J_0 \leq nv(n)\) discovered by Ore \cite{Ore}.

We define \(P\) to be \emph{weakly valid}\footnote{This is new.} if it satisfies these same conditions but with \(s\) restricted to \(\{s_t\}\), i.e.
\begin{align*}
& p^{s_t} \leq b_t \\
b_t = n \implies& \ell_P(n,s_t) = 0 &\text{(Ore 1)} \\
& \ell_P(n,s_t) \leq 0 &\text{(Ore 2)} \\
b_t < n \implies& \ell_P(b_t,s_t) \geq 1 &\text{(Ore 3)}\\
b_r = b_t < n \implies& \ell_P(b_t,s_t) = \ell_P(b_r,s_r) &\text{(Consistency)} \\
p^{s_r} \leq b_t < n \implies& \ell_P(b_t,s_t) \geq \ell_P(b_t,s_r). &\text{(Bounding)} \\
\end{align*}

Noting that \(\ell_P(i,s_t) = a_t + \ceiling{\tfrac{b_t-i}{n}} - B(i,p^{s_t}) + 1\) is actually a function of \(i\) and the vertex \((p^{s_t}, J_t=a_tn+b_t)\), we deduce that if the vertices of \(P'\) are a subset of the vertices of \(P\) and \(P\) is weakly valid, then \(P'\) is also weakly valid. In particular if \(P'\) is valid, then removing vertices gives a weakly valid polygon.

\subsection{Enumeration}
\label{ext-sec-rampol-enum}

We can use the Ore bounds to enumerate all possible ramification polygons \(P\) for a given degree \(n\). We perform a branching algorithm to assign all possible combinations of vertices. The state of a branch is a paritally-assigned polygon \(P\) and an integer \(0 \leq S \leq v_p(n)\). First, we branch over each \(J_0=0,1,\ldots,n v(n)\) satisfying the Ore bound, and set the state to \(P=[(1,J_0),(p^{v_p(n)},0),(n,0)]\) and \(S=1\). Given the state \(P=[(p^{s_0}=1,J_0),\ldots,(p^{s_k},J_k),(p^{v_p(n)},0),(n,0)]\) and \(s_k<S<v_p(n)\) we consider adding a vertex of the form \((p^S,*)\) to \(P\) or not. Hence we branch: either we don't add a vertex, in which case the state of the branch becomes \(P\) and \(S+1\); otherwise we branch for each possible new vertex \((p^S,J)\) with \(0 < J < P(p^S)\), in which case the state becomes \(P+(p^s,J)\) and \(S+1\). Finally, if we are given the state \(P\) and \(S=v_p(n)\) then we have decided on all vertices for \(P\), and so we check it is valid and if so, output it.

This algorithm is quite impractical but can be made practical by terminating a branch if the current \(P\) is not weakly valid. We know that removing vertices from a valid polygon preserves weak validity, so if \(P\) is not weakly valid it can not be augmented to a valid one.\footnote{This algorithm is similar to one described in Sinclair's thesis \cite[\S3.3]{SinclairTh}. It is not clear if that algorithm is correct because it does not seem to use our concept of weak validity, which appears necessary for building up polygons one vertex at a time. Furthermore, tables \cite{SinclairTables} produced from their implementation are missing some polygons. For example, the tables for degree 8 extensions of \(\QQ_2\) are missing the polygon \([(1,7),(2,6),(4,4),(8,0)]\), which is the ramification polygon of \(x^8+2x^7+2x^6+2x^4+2\).}

For example, to compute all 447 ramification polygons of degree 16 over \(\QQ_2\) would require considering around 300,000 branches, but can reduce this to 1602 branches by terminating using weak validity. The 6849 ramification polygons of degree 32 are found with 29,730 branches instead of around 600,000,000.

\subsection{Template}
\label{ext-sec-rampol-template}

Suppose we write \(f_i = \sum_{0 \leq k \leq \infty} \hat f_{i,k} \pi^k\) for \(f_{i,k} \in \FF_K\), where \(\hat\cdot:\FF_K \to \cO_K\) is a choice of representative. A \emph{template} is a collection of sets \(X_{i,k} \subset \FF_K\) defining a set of monic polynomials\footnote{Enumerating Eisenstein polynomials via templates is used throughout \cite{PS}, and in particular in Algorithm 6.1.}
\[X = \braces{ x^n + \sum_{i=0}^{n-1} x^i \sum_{k=0}^\infty \hat f_{i,k} \pi^k \,:\, f_{i,k} \in X_{i,k} } \subset K[x].\]

For example, the template for all Eisenstein polynomials has \(X_{i,0} = \{0\}\), \(X_{1,1} = \FF_K^\times\), otherwise \(X_{i,k} = \FF_K\).

Given a valid ramification polygon \(P\), a template consisting precisely of the Eisenstein polynomials with ramification polygon \(P\) is given by the template for Eisenstein polynomials with additionally \(X_{i,k} = \{0\}\) for all \(p^s \leq i\), \(k < \ell_P(i,s)\), and \(X_{b_t, \ell_P(b_t,s_t)} = \FF_K^\times\) for all \(t\).

Given a finite template, one can enumerate all of its polynomials by enumerating the cartesian product of the \(X_{i,k}\).

One can show as a consequence of Krasner's lemma that if \(f,f'\in K[x]\) are Eisenstein and \(f-f'\) has coefficients of valuation more than \(1+2J_0/n\) then they generate the same field. Therefore a template can be made finite by setting \(X_{i,k} = \{0\}\) for \(k>1+2J_0/n\), and its polynomials will between them generate the same set of extensions. In \S\ref{ext-sec-shrink} we shall see how to make the template even smaller.

\section{Fine ramification polygon}
\label{ext-sec-finerampol}

\subsection{Definition}
\label{ext-sec-finerampol-def}

We define the set \(\cPst = \cP \cap P\) to be the \emph{fine ramification polygon of \(f\)}. It is not strictly a polygon itself, but its lower convex hull is \(P\), and so this is a finer quantity than the ramification polygon because it can include points in the interiors of the faces. As discussed, all points \((j,*)\) for \(j \leq p^{v_p(n)}\) in \(\cPst\) are of the form \((p^s,*)\) and so we denote \(\cPst\) by
\[\cPst = [(p^{s_0}=1,J_0),\ldots,(p^{s_u},J_u=0),\ldots,(n,0)].\]
Note that the quantities \(u\), \(s_t\) and \(J_t\) may be different to their earlier meaning in the context of the plain ramification polygon, since there may now be extra points between the vertices. Also note that there are possibly some extra points between \((p^{s_u},0)\) and \((n,0)\).

\subsection{Invariant}
\label{ext-sec-finerampol-inv}

We shall see in \S\ref{ext-sec-residues-inv} that \(\cPst\) is an invariant of \(L/K\).

\subsection{Validity}
\label{ext-sec-finerampol-valid}

A potential fine ramification polygon \(\cPst\) is valid if and only if there are \(F_i \in \ZZ\) such that \(R_j = J\) for all points \((j,J) \in \cPst\) and for all \(j\) without a point above \(j\) we have \(R_j > \cPst(j)\).

First we consider the horizontal face. For \(p^{s_u} \leq j \leq n\) we need that \(R_j = 0\) if and only if \((j,0) \in \cPst\). Observe that \(r_j \equiv \binom{n}{j}\) and therefore \(R_j = 0\) if and only if \(B(n,j) = 0\). Hence the condition for this range of \(j\) is that \((j,0) \in \cPst\) if and only if \(B(n,j)=0\).

For the remaining points \((p^s,*)\), the analysis is almost identical to the case with the plain ramification polygon, the only difference being that we require \(R_{p^s} > \cPst(p^s)\) whenever there is not a point \((p^s,*) \in \cPst\). Following the same steps, we deduce the following conditions:
\begin{align*}
& F_0 = 1 \\
0 < i < n \implies& F_i \geq 1 \\
& F_n = 0 \\
p^s \leq i \leq n \implies& F_i \geq \ell_{\cPst}(i,s) \\
& F_{b_t} = \ell_{\cPst}(b_t,s_t) \\
\end{align*}
where
\[\ell_{\cPst}(i,s_t) := a_t - B(i,p^{s_t}) + 1 + 1[i<b_t]\]
and for \((p^s,*) \not\in \cPst\)
\[\ell_{\cPst}(i,s) := \floor{\frac{P(p^s)-i}{n}} - B(i,p^s) + 2.\]

Observing that these conditions are essentially identical to those from the previous section, up to the redefinition of \(\ell_{\cPst}\), then we similarly find that \(\cPst\) is valid if and only if\footnote{This is essentially \cite[Prop. 3.9]{PS}, which is slightly mis-stated because it does not include the Ore 2 condition in the case that there is no point \((p^s,*) \in \cPst\).}
\begin{align*}
B(n,j)=0 \iff& (j,0) \in \cPst &\text{(Tame)}\\
& p^{s_t} \leq b_t \\
b_t = n \implies& \ell_{\cPst}(n,s_t) = 0 &\text{(Ore 1)} \\
& \ell_{\cPst}(n,s) \leq 0 &\text{(Ore 2)} \\
b_t < n \implies& \ell_{\cPst}(b_t,s_t) \geq 1 &\text{(Ore 3)} \\
b_r = b_t < n \implies& \ell_{\cPst}(b_t,s_t) = \ell_{\cPst}(b_r,s_r) &\text{(Consistency)} \\
p^s \leq b_t < n \implies& \ell_{\cPst}(b_t,s_t) \geq \ell_{\cPst}(b_t,s). &\text{(Bounding)} \\
\end{align*}

If \(\cPst\) is valid, then \(f(x)\) has \(\cPst\) as its fine ramification polynomial if and only if \(F_{b_t} = \ell_{\cPst}(b_t,s_t)\) for all \(t\) and \(F_i \geq \ell_{\cPst}(i,s)\) whenever \(p^s \leq i \leq n\).\footnote{This is essentially \cite[Prop. 3.10]{PS}.}

We define \(\cPst\) to be \emph{weakly valid} if these conditions hold only for \(s \in \{s_t\}\). As before, weak validity is preserved under removing points from \(\cPst\).

Enumerating the possible fine ramification polygons and producing their templates is essentially the same as for plain ramification polygons.

\section{Residues}
\label{ext-sec-residues}

\subsection{Definition}
\label{ext-sec-residues-def}

We define \(\rho_j\) to be the leading \(\alpha\)-adic coefficient of \(r_j\). Recall that if \(R_j := n v(r_j) = an+b\) then \(r_j \sim \binom{b}{j} f_b \alpha^{b-n}\). Hence
\[\rho_j \equiv r_j \alpha^{-R_j} \equiv \binom{b}{j} (\phi_b \pi^{F_b}) (-\phi_0 \pi)^{-1-a}\]
where \(\phi_i \equiv f_i \pi^{-F_i}\).

We refer to \(-\phi_0\) as the \emph{uniformizer residue of \(\alpha\) (with respect to \(\pi\))} because \(\alpha^n \sim -f_0 \sim (-\phi_0) \pi\).

We extend the points of the fine ramification polygon to include the residues \((j,R_j,\rho_j)\), giving the \emph{fine ramification polygon with residues}\footnote{In \cite{PS}, the equivalent notion of \emph{residual polynomial} is studied instead, but this additional structure is not actually used.}
\[\cPstres = [(p^{s_0}=1, J_0, \gamma_0),\ldots,(p^{s_u}, 0, \gamma_u), \ldots, (n,0,1)].\]

\subsection{Invariant}
\label{ext-sec-residues-inv}

If the ramification polygon \(P\) has a face \(((j_0,R_{j_0}),(j_1,R_{j_1}))\) of width \(w=j_1-j_0\) and slope \((R_{j_1}-R_{j_0})/w=-h/e\), we define its \emph{residual polynomial} to be
\[A(x) \equiv r_{j_1}^{-1} \sum_{i=0}^{w/e} r_{j_0+ie} \alpha^{jh-R_{j_0}} x^i \equiv \sum_{\substack{(j,R_j,\rho_j) \in \cPstres \\ j_0 \leq j \leq j_1}} \rho_j x^{(j-j_0)/e} \in \FF_K[x].\]

It is well-known that the \(w\) roots \(\tfrac{\alpha-\sigma\alpha}{\alpha}\) of \(r\) of valuation \(h/e\) satisfy
\[A((\tfrac{\alpha-\sigma\alpha}{\alpha})^e / \alpha^h) \equiv 0.\]

Observe that \(A\), being monic, is determined by its roots. Recalling that if we choose another uniformizer of \(\alpha' \in L=K(\alpha)\) such that \(\alpha' \sim \delta \alpha\) then \(\tfrac{\alpha'-\sigma\alpha'}{\alpha'} \sim \tfrac{\alpha-\sigma\alpha}{\alpha}\), these are the roots of \(r(x)\) and \(r'(x)\), and we deduce that if \(A'\) is the corresponding residual polynomial then \(A'(x) \equiv A(\delta^h x)\) where \(\alpha' \sim \delta \alpha\).

Firstly this implies that a coefficient of \(A'\) is zero if and only if the corresponding coefficient of \(A\) is zero. This implies that \(\cPst=\cPstpr\) is an invariant of \(L/K\).

Furthermore, this implies that \(\rho'_j \equiv \rho_j \delta^{-R_j}\) for \((j,R_j,\rho_j) \in \cPstres\). Therefore, given \(\cPstres\) and \(\cPstrespr\), we consider them equivalent if they are equal as fine ramification polygons and if there exists \(\delta \in \FF_K^\times\) such that \(\rho'_j \equiv \rho_j \delta^{-R_j}\). Then equivalence classes are an invariant of \(L/K\).\footnote{This is equivalent to the invariant \(\cA\) of \cite[\S4.1]{PS}.}

Note that \(\rho'_j \equiv \rho_j \delta^{-R_j}\) if and only if \(r'_j \sim r_j\) when \(\alpha' \sim \delta \alpha\). Therefore being equivalent is the same thing as \(v(r_j - r'_j) > P(j)\) for all \(1 \leq j \leq n\).

\subsection{Validity}
\label{ext-sec-residues-valid}

Given a fine ramification polygon with residues \(\cPstres\), then it is valid if and only if: (a) it is valid as a fine ramification polygon; (b) for \((j,0,\rho_j) \in \cPstres\) we have \(\rho_j \equiv r_j \equiv \binom{n}{j}\); and (c) there are \(\phi_i \in \FF_K^\times\) such that for each \(t\)
\[\gamma_t \equiv \binom{b_t}{p^{s_t}} (\phi_{b_t} \pi^{F_{b_t}}) (-\phi_0 \pi)^{-1-a_t} \equiv \beta(b_t,p^{s_t}) \phi_{b_t} (-\phi_0)^{-1-a_t}\]
where \(\beta(i,j) \equiv \binom{i}{j} \pi^{-B(i,j)}.\)

Eliminating \(\phi_i\) from these latter equations for \(0<i \leq n\), we find they are soluble if and only if there is \(\phi_0 \in \FF_K^\times\) such that\footnote{This is essentially \cite[Prop. 4.5]{PS}, which is slightly mis-stated because it misses the conditions in the case \(b_t=n\).}
\begin{align*}
b_t = n \implies& \gamma_t \equiv \beta(n,p^{s_t}) (-\phi_0)^{-a_t-1} \\
b_t = b_r < n \implies& \frac{\gamma_t}{\gamma_r} \equiv \frac{\beta(b_t,p^{s_t})}{\beta(b_r,p^{s_r})} (-\phi_0)^{a_r-a_t}.
\end{align*}

If \(\cPstres\) is valid, then \(f(x)\) has \(\cPstres\) as its fine ramification polygon with residues if and only if \(\phi_0\) is a solution to these equations and
\[\phi_{b_t} \sim \gamma_t \beta(b_t,p^{s_t})^{-1} (-\phi_0)^{a_t+1}\]
for each \(t\).

\subsection{Enumeration}
\label{ext-sec-residues-enum}

Given \(\cPst\), we can enumerate all possible \(\cPstres\) extending it by initially having \(\gamma_t\) unset for all \(t\). Similar to our branching algorithm for enumerating ramification polygons, we branch for each \(t\) over the possible values of \(\gamma_t\).

We can make this practical by terminating a branch if the current partial assignment of \(\gamma_t\) is not \emph{weakly valid}, where a partial assignment is weakly valid if there exists \(\phi_0 \in \FF_K^\times\) such that the above conditions involving only the assigned \(\gamma_t\) are true.

Suppose we have assigned \(\gamma_0,\ldots,\gamma_{t-1}\) so far and are branching over possibilities for \(\gamma_t\). We can restrict the algorithm to produce one representative \(\cPstres\) per class by considering \(\gamma_t\) equivalent to \(\gamma'_t\) if there exists \(\delta \in \FF_K^\times\) such that \(\delta^{-J_k} \equiv 1\) for \(k < t\) and \(\gamma'_t \equiv \gamma_t \delta^{-J_t}\), and only branching over representatives of equivalence classes.

\subsection{Template}
\label{ext-sec-residues-template}

The choice of \(\phi_0\) means there is not a template (by our definition) for Eisenstein polynomials corresponding to \(\cPstres\). This will be fixed in the next section by including \(\phi_0\) in the invariant.

\section{Uniformizer residue}
\label{ext-sec-constcoeff}

\subsection{Invariant}
\label{ext-sec-constcoeff-inv}

We now consider the extent to which the pair \((\cPstres, \phi_0)\) is an invariant of \(L/K\). We have seen that changing uniformizer \(\alpha \to \alpha' \sim \delta \alpha\) changes \(\rho_j \to \rho'_j \sim \delta^{-R_j} \rho_j\). Since by definition \(\alpha^n \sim -\phi_0 \pi\), we find that \(\phi_0 \to \phi'_0 \sim \delta^n \phi_0\).

Therefore we say \((\cPstres, \phi_0)\) and \((\cPstrespr, \phi'_0)\) are equivalent if there is \(\delta \in \FF_K^\times\) such that \(\rho'_j \sim \delta^{-R_j} \rho_j\) for all \((j,R_j,\rho_j) \in \cPstres\) and \(\phi'_0 \sim \delta^n \phi_0\). Equivalence classes are an invariant for \(L/K\).\footnote{This is equivalent to the invariant \(\cA^*\) of \cite[\S4.2]{PS}.}

By the preceding section, \(\phi_{b_t}\) is determined for all \(t\) by the choice of representative. Therefore considering \((\cPstres, \phi_0, \{\phi_{b_t}\,:\,t\})\) does not give a finer invariant.

\subsection{Validity and enumeration}
\label{ext-sec-constcoeff-val}
\label{ext-sec-constcoeff-enum}

Suppose we are given a valid \(\cPstres\), then the valid \(\phi_0\) which extend this are precisely the solutions to the system of the previous section. Hence to enumerate them all, we find the solutions.


If we want to find representatives of equivalence classes, note that \((\cPstres, \phi_0)\) and \((\cPstres, \phi'_0)\) are equivalent if there is \(\delta\) such that \(\delta^{-g} = 1\) and \(\delta^n = \phi'_0 / \phi_0\) where \(g = \gcd_t J_t = \gcd_{(j,R_j) \in \cPst} R_j\).
We use this criterion to test for equivalence, and so choose one \(\phi_0\) per equivalence class.

\subsection{Template}
\label{ext-sec-constcoeff-template}

The template for Eisenstein polynomials corresponding to \((\cPstres, \phi_0)\) is the template for \(\cPst\) with the changes \(X_{0,1} = \{\phi_0\}\) and
\[X_{b_t,\ell_{\cPst}(b_t,s_t)} = \{\gamma_t \beta(b_t,p^{s_t})^{-1} (-\phi_0)^{a_t+1}\}.\]

\section{Change of uniformizer}
\label{ext-sec-shrink}

So far we have studied the change of uniformizer \(\alpha \to \alpha' = \delta \alpha'\) and the effect of \(\bar\delta\). We now consider smaller perturbations, i.e. \(\alpha' = \alpha(1 + u \alpha^m)\) for \(m>0\), \(u \in \cO_L\). These will allow us to considerably reduce the size of our templates.

Letting \(f'(x)\) be the minimal polynomial for \(\alpha'\) then
\[f(\alpha')-f(\alpha') = f(\alpha') = \alpha^n r(\tfrac{\alpha'}{\alpha}-1) = \alpha^n r(u \alpha^m)\]
and
\[f(\alpha') - f'(\alpha') = \sum_{i=0}^{n-1} (f_i - f'i)(\alpha')^i.\]

Writing
\[S_m(u) \equiv \alpha^{-C_m} r(\alpha^m x) \in \FF_K[x]\]
with \(C_m\) as large as possible, and writing \(C_m = c_m n + d_m\) with \(0 \leq d_m < n\) then
\[S_m(u) \equiv (f_{d_m} - f'_{d_m}) \alpha^{d_m - n - C_m} \equiv (f_{d_m} - f'_{d_m}) (-\phi_0 \pi)^{- 1 - c_m}.\]

For \(m>0\), \(S_m\) only consists of terms from points on faces of \(\cPst\) with negative slope, i.e. the points \((p^{s_t},J_t)\), and hence \(S_m\) only has terms at powers of \(p\) and so is additive. Therefore it defines a \(\FF_p\)-linear map \(S_m : \FF_K \to \FF_K\).

By changing uniformizer \(\alpha \to \alpha'\) we change \(f_{d_m,1+c_m}\) by \((-\phi_0)^{1+c_m} S_m(u)\). Therefore if we restrict the template by setting \(X_{d_m,1+c_m}\) to be a set of coset representatives of \(\FF_K^+ / (-\phi_0)^{1+c_m} S_m(\FF_K^+)\), then it will produce the same set of fields.

Observe that since \(C_m\) is strictly increasing in \(m\) (in fact \(C_m\) is essentially the Hasse-Herbrand transition function; see e.g. \cite{SerLF,Helou}) then this change to the template is independent to the changes made for any higher \(m\), and so we can make all of these changes independently.

In particular, if \(m\) is higher than the negative of the slope of the steepest face of \(P\) then \(S_m(x) = x\) is surjective, and so we may set \(X_{d_m,1+c_m} = \{0\}\). In this region \(C_m = J_0 + m\), and so this sets all but finitely many \(X_{i,k} = \{0\}\). Note that this change is independent of \(\phi_0\) (unlike for smaller \(m\)) and therefore can be used when enumerating polynomials from simpler invariants like \(P\) or \(\cPst\).

\section{Implementation notes}
\label{ext-sec-impl}

\subsection{Representation of invariants}

Our representation of a ramification polygon is a little more general than is presented in this paper. We actually represent a polygon as a list
\[[(x_0,J_0,\sim_0,\rho_{x_0}),(x_1,J_1,\sim_1,\rho_{x_1}),\ldots,(n,0,*,\rho_n)]\]
where \(x_i = p^i\) for \(i \leq v_p(n)\). That is, we have a point for every power of \(p\). The relation \(\sim_i\) specifies that \(R_{x_i} \sim_i J_i\) where \(\sim_i\) is \(=\), \(\geq\) or \(>\). Hence if \(\sim_i\) is \(=\) then \((x_i,J_i)\) is a point of the (fine) ramification polygon; if it is \(>\) then it is not; if it is \(\geq\) then it is unspecified.

In particular, for the plain ramification polygon, \(\sim_i\) is \(=\) at the vertices and \(\geq\) elsewehere. For the fine ramification polygon, \(\sim_i\) is \(=\) or \(>\).

The residues \(\rho_{x_i}\) may be unspecified. Where they are specified, they must be non-zero and \(\sim_i\) must be \(=\).

The arguments from earlier sections are simply modified to yield validity and equivalence conditions for this more general object. One way in which this is useful is that we can now specify one residue at a time, and check for validity at each step, and therefore the branching algorithms for enumerating these invariants may be implemented.

\subsection{Consistency of roots}

Checking for equivalence of ramification polygons with residues requires solving a system of equations of the form \(x^{k_i} = a_i\).

We compute the extended GCD
\[K := \gcd(k_1,k_2,\ldots) = \sum_i b_i k_i\]
so that
\[x^K = \prod_i x^{b_i k_i} = \prod_i a_i^{b_i} =: A\]
and
\[a_i = x^{k_i} = x^{K(k_i/K)} = A^{k_i/K}.\]

We deduce the system is solvable if and only if \(a_i = A^{k_i/K}\) and \(x^K = A\) is solvable. Hence we can check for consistency and reduce the system to a single equation.

To limit the solutions to \(\FF_K^\times\), we let \(q=\abs{\FF_K}\) and include \(x^{q-1} = 1\) in the system of equations.

\subsection{Binomial coefficients}

A little care is needed to compute with binomial coefficients \(\pi\)-adically.

One can see that
\[v_p(k!) = \sum_{i=1}^{\infty} \floor{\tfrac{k}{p^i}}\]
and from which we can compute
\[B(r,k) = e(K/\QQ_p)(v_p(r!) - v_p(k!) - v_p((r-k)!)).\]

By Wilson's theorem, if \(k = ap+b\) then the ``\(p\)-unit'' part of the factorial is
\[U_p(k) := \prod_{1 \leq i \leq k, v_p(i)=0} i \equiv (-1)^a b! \mod p\]
from which we can compute the ``\(p\)-shifted factorial''
\[S_p(k) := k! / p^{v_p(k!)} \equiv \prod_{i=1}^{\infty} U_p(\floor{\tfrac{k}{p^i}}) \mod p\]
and hence
\[\beta(r,k) \equiv \frac{S_p(r)}{S_p(k) S_p(r-k)} \gamma^{-v_p \binom{r}{k}}\]
where \(\pi^{e(K/\QQ_p)} \sim \gamma p\) (i.e. \(\gamma\) is the uniformizer residue of \(\pi\) with respect to \(p\) c.f. \S\ref{ext-sec-residues-def}).

\bibliography{refs}

\end{document}